\documentclass[reqno]{amsart}
\usepackage{hyperref, amssymb, mathrsfs,dsfont}

\begin{document}

\vskip -8in
\begin{center}
	\rule{8.5cm}{0.5pt}\\[-0.1cm] 
	{\small To appear in 
		{\it Trends in Mathematics},  Birkh\"auser,
		2024.
	}
	%%%%%%%%%%%%%%%%%
	%{\small Preprint
		%}%ALERT%%
	%%%%%%%%%%%%%%%%%
	\\[-0.25cm] \rule{8.5cm}{0.5pt}
\end{center}
\vspace {2.2cm}

\title[Critical Sobolev embeddings in the Heisenberg group]
{Global Compactness, subcritical approximation of the Sobolev quotient, and a related concentration result in the Heisenberg group}

\author[Giampiero Palatucci, Mirco Piccinini, Letizia Temperini\hfil \hfilneg] {Giampiero Palatucci, Mirco Piccinini, Letizia Temperini}  

\address{Giampiero Palatucci, Mirco Piccinini, Letizia Temperini\newline
Dipartimento di Scienze Matematiche, Fisiche e Informatiche, Universit\`a di Parma\\ Campus - Parco Area delle Scienze 53/a, 43124 Parma, Italy}
\email{giampiero.palatucci@unipr.it,  mirco.piccinini@unipr.it, letizia.temperini@outlook.it}

\begin{abstract}
We investigate some effects of the lack of compactness in the critical Sobolev embedding in the Heisenberg group.\end{abstract}

\maketitle \numberwithin{equation}{section}
\newtheorem{theorem}{Theorem}[section]
\newtheorem{corollary}[theorem]{Corollary}
\newtheorem{lemma}[theorem]{Lemma}
\newtheorem{remark}[theorem]{Remark}
\newtheorem{problem}[theorem]{Problem}
\newtheorem{example}[theorem]{Example}
\newtheorem{definition}[theorem]{Definition}
\allowdisplaybreaks

%%%% CUSTOM DEFINITIONS
\def \eps{{\varepsilon}}
\def\Om{\Omega}
\def \h{{\mathds{H}^n}}
\def \r{{\mathds{R}}}
\def \N{{\mathds{N}}}
\newcommand{\El}{\mathcal{E}_\lambda}
\newcommand{\Es}{\mathcal{E}^*}
\newcommand{\Sc}{{S}^1_0}
\newcommand{\Sob}{S^\ast}
\newcommand{\Som}{\Sob_\Omega}
\def\mea{\mathcal{M}(\overline{\Omega})}
\newcommand{\snr}[1]{\lvert #1\rvert}
\newcommand{\Car}{\Sigma_{\Omega,D_H}}
\newcommand{\Xc}{{X}}
\newcommand{\Omb}{\overline{\Omega}}
\newcommand{\Irm}{\text{I}}
\newcommand{\towt}{\stackrel{\mathcal{T}}{\rightarrow}}
\newcommand{\tows}{\stackrel{\ast}{\rightharpoonup}}
\newcommand{\Fep}{\mathcal{F}_{\varepsilon}}
\def\dys{\displaystyle}
\newcommand{\Fc}{\mathcal{F}}
\def\gamp{\Gamma^{+}}
\def\ue{u_{\eps}}
\def\mue{\mu_{\eps}}
\newcommand{\uz}{u^{(0)}}
\newcommand{\lambdakj}{{\lambda_k^{(j)}}}
\newcommand{\xikj}{{\xi_k^{(j)}}}
\newcommand{\ukj}{{u_k^{(j)}}}
\newcommand{\rkj}{{r_k^{(j)}}}
\newcommand{\uj}{{u^{(j)}}}
\def\mea{\mathcal{M}(\overline{\Omega})}
\newcommand{\twe}{T{w}_{\eps}}
\newcommand{\uer}{u_{\eps,r}}
\newcommand{\tuer}{T{u}_{\eps,r}}
\newcommand{\Ssub}{\Sob_\eps}
%%%%%%%%%%%

\vspace{-5mm}  
\section{Critical Sobolev embeddings in the Heisenberg group}\label{sec_cri}
Let  $\h := (\mathds{C}^n\times\mathds{R},\circ, \delta_\lambda)$ be the usual Heisenberg-Weyl group, endowed with the group multiplication law~$\circ$,
 $$
   	\xi \circ \xi' := \Big(x+x',\, y+y',\, t+t'+2\langle y,x'\rangle-2\langle x,y'\rangle \Big) \ \,
	$$
	for $\xi:=(x+iy, t)$  and $\xi':=(x'+iy', t') \in \r^n\times \r^n \times \r$,
whose group of non-isotropic {\it dilations}~$\{\delta_\lambda\}_{\lambda>0}$ on~$\r^{2n+1}$ is given by
 \begin{equation}\label{def_philambda}
   	                 \xi   \mapsto \delta_\lambda(\xi):=(\lambda x,\, \lambda y,\, \lambda^2 t).
	               \end{equation}

Consider the standard Folland-Stein-Sobolev space~$\Sc(\h)$ defined as the completion of~$C^\infty_0(\h)$ with respect to the homogeneous subgradient norm~$\|D_H\cdot\|_{L^2}$, where the  horizontal (or intrinsic) gradient~$D_H$ is given by
   $$
   D_H u(\xi) := \big( Z_1 u(\xi),\dots, Z_{2n}u(\xi)\big),
   $$
   with  $Z_j := \partial_{x_j} +2y_j\partial_t$,  $Z_{n+j}:= \partial_{y_j}-2x_j\partial_t$ for $1 \leq j\leq n$,  and   $T:=\partial_t$ being the Jacobian base of the Heisenberg Lie algebra.

\vspace{1mm}

As well known, the following Sobolev-type inequality  holds for some positive constant~$\Sob$,
\begin{equation}\label{folland}
\|u\|^{2^\ast}_{L^{2^\ast}} \leq \Sob \|D_Hu\|^{2^\ast}_{L^2}\,, \quad \forall u \in S^1_0(\h)\,,
\end{equation}
where  $2^\ast=2^\ast(Q):=2Q/(Q-2)$ is the Folland-Stein-Sobolev critical exponent, depending on the {\it homogeneous dimension}~$Q:=2n+2$ of the Heisenberg group~$\h$.

The validity of~\eqref{folland} is equivalent to show that the constant~$\Sob$ defined in the 
following maximization problem,
\begin{equation}\label{critica0}
\Sob:=\sup\left\{\int_{\h}|u(\xi)|^{2^\ast}\,{\rm d}\xi \, : \, u\in \Sc(\h), \int_{\h}|D_H u(\xi)|^2{\rm d}\xi \leq 1\right\},
\end{equation}
is finite.
The explicit form of the maximizers has been showed, amongst other results, in the breakthrough paper by Jerison and Lee~\cite{JL88}, together with the computation of the optimal constant in~\eqref{critica0}. 

For any bounded domain~$\Omega\subset\h$, consider now
\begin{equation}\label{critica}
\Som:=\sup\left\{\int_{\Omega}|u(\xi)|^{2^\ast}\,{\rm d}\xi \, : \, u\in \Sc(\Omega), \int_{\Omega}|D_H u(\xi)|^2{\rm d}\xi \leq 1\right\},
\end{equation}
where the Folland-Stein-Sobolev space $\Sc(\Omega)$ is given by the closure of $C^\infty_0(\Omega)$ with respect to the homogeneous subgradient norm in $\Omega$. One can check that $\Som\equiv \Sob$ via a standard
scaling argument, and thus -- in view of the explicit form of the optimal functions in~\eqref{critica0} --  the variational problem~\eqref{critica} has no maximizers. The situation changes drastically for the subcritical embeddings:
$\Sc(\Omega) \hookrightarrow L^{2^\ast-\eps}(\Omega)$ is compact (for each~$0<\eps<2^\ast-2$), and this  guarantees the existence of a maximizer~$u_\eps\in \Sc(\Omega)$ for
   \begin{equation}\label{sobolev}
\Sob_\eps:=\sup\left\{\int_{\Omega}|u(\xi)|^{2^\ast-\eps}\,{\rm d}\xi \, : \, u\in \Sc(\Omega), \int_{\Omega}|D_H u(\xi)|^2{\rm d}\xi \leq 1\right\}.
\end{equation}
Such a dichotomy can be also found in the Euler-Lagrange equation for the energy functionals in~\eqref{sobolev}; that is,
\begin{equation}\label{equazione}
-\Delta_H u_\eps = \lambda |u_\eps|^{2^\ast-\eps-2} u_\eps \, \ \text{in} \ (\Sc(\Omega))',
\end{equation}
where $\lambda$ is a Lagrange multiplier, and $  \Delta_H :=  \sum_{j=1}^{2n} Z^2_j$ is the standard Kohn Laplacian (or sub-Laplacian) operator. While when $\eps>0$ it has a solution~$u_\eps$, the problem above becomes very delicate when $\eps=0$: one falls in the CR~Yamabe equation realm, and even the existence of the solutions is not granted. 
In view of such a qualitative change when $\eps=0$ \big(in both~\eqref{sobolev} and \eqref{equazione}\big), it sounds natural to analyze the asymptotic behavior as $\eps$ goes to $0$ of both
the subcritical Sobolev constant~$\Sob_\eps$ in the Heisenberg group given in~\eqref{equazione} and of the corresponding optimal functions~$u_\eps$ of the embedding $\Sc(\Omega) \hookrightarrow L^{2^\ast-\eps}(\Omega)$. This is the aim of the papers~\cite{PP23h} and~\cite{PPT23}, whose main results will be stated in the rest of the present note.

\vspace{1mm}
\section{Subcritical approximation of the Sobolev quotient}\label{sec_gamma}
 Our first  result is the subcritical approximation of the Sobolev embedding~$S^\ast$ in the Heisenberg group  described 
  below. 
    \begin{theorem}[See Theorem~1.1 in~\cite{PPT23}]\label{the_gamma-intro} 
Let $\Omega\subseteq\h$ be a bounded domain, and denote by~$\mea$ the set of nonnegative Radon 
measures in~$\Omega$. Let $\Xc=\Xc(\Omega)$ be the space
$$
\Xc:=\Big\{(u,\mu) \in \Sc(\Omega)\times\mea: \mu \geq |D_H u|^{2}{\rm d}\xi, \, \mu(\Omb)\leq 1\Big\},
$$ endowed with the product  topology $\mathcal{T}$ such that
\begin{equation}\label{def_top1}
(u_k,\mu_{k}) \towt (u,\mu) \ \, {\stackrel{\text{def}}{\Leftrightarrow}} \ \, \begin{cases} u_{k} \rightharpoonup u \ \text{in} \ L^{2^{\ast}}\!(\Omega), \\
\mu_{k} \tows \mu \ \text{in} \ \mea.
\end{cases}
\end{equation}
\noindent
\\ Let us consider the following family of functionals,
\begin{equation*}\label{def_fue1}
\Fep(u,\mu):= \dys \int_{\Omega}|u|^{2^{\ast}\!-\eps} {\rm d}\xi  \ \ \ \forall (u,\mu) \in \Xc\,.
\end{equation*}
Then, as $\eps\to 0$, the $\Gamma^{+}$-limit of the family of functionals $\Fep$ with respect to the topology~$\mathcal{T}$ given by~{\rm(\ref{def_top1})} is the functional $\Fc$ defined by
\begin{equation*}\label{def_fu}
\Fc(u,\mu)=\int_{\Omega}|u|^{2^{\ast}}{\rm d}\xi + S^{\ast}\sum_{j=1}^{\infty}\mu_{j}^{\frac{2^{\ast}}{2}} \ \ \ \forall (u,\mu) \in \Xc.
\end{equation*}
Here $S^{\ast}$ is the best Sobolev constant in $\h$, $2^\ast=2Q/(Q-2)$ is the Folland-Stein-Sobolev critical exponent, 
and the numbers $\mu_j$ are the coefficients of the atomic part of the measure $\mu$.
\end{theorem}
  
In order to prove such a result in the very general situation considered here, and thus requiring no additional regularity assumptions nor special geometric features on the domains,  we attack the problem pursuing a new approach and for this we rely on De Giorgi's $\Gamma$-convergence techniques. This is in the same spirit of previous results regarding the classical Sobolev embedding in the Euclidean framework, as seen in~\cite{AG03,Pal11,Pal11b}, 
  though the core of the proof in~\cite{PPT23} goes in a very different line because the optimal recovery sequences have been concretely constructed whereas in all the aforementioned Euclidean papers such an existence result has been proven via compactness and locality properties of the $\Gamma$-limit energy functional. In this respect, the adopted strategy is surprisingly close to that in the  
  {\it fractional Sobolev spaces\,} framework (~\!\!\!\!\:\cite{
   PP14,
   PPS15}), but various differences evidently arose because of the natural discrepancy
between the involved frameworks.

It could be interesting to investigate whether or not the techniques introduced in~\cite{PPT23} and~\cite{PPS15} could be combined with the estimates involving the ``nonlocal tail'' in the Heisenberg framework firstly introduced in~\cite{PP22} in order to prove a similar result for fractional Folland-Stein-Sobolev spaces; see also~\cite{MPPP23,Pic22,PP23}.
\vspace{2mm}
  
    As a corollary of Theorem~\ref{the_gamma-intro}, one can deduce that the sequences of maximizers~$\{u_\eps\}$ for the subcritical Sobolev quotient~$\Sob_\eps$  concentrates energy at one point~$\xi_{\rm o}\in\Omb$, and this is in clear accordance 
  with the analogous result in the Euclidean case.

   \vspace{2mm}   
  \begin{theorem}[See Theorem~1.2 in~\cite{PPT23}]\label{cor_concentration}   	Let $\Om \subset \h$ be a bounded domain and 
  let~$\ue\in \Sc(\Omega)$ be a maximizer for~$S^{\ast}_{\eps}$.  Then, as~$\varepsilon=\varepsilon_k \to 0$, up to subsequences, we have that
  	there exists $\xi_{\rm o} \in \Omb$ such that
  	$$
  	u_k=u_{\eps_k}  \rightharpoonup 0  \ \mbox{in} \ L^{2^\ast}\!(\Omega),
  	$$
  	and 
  	$$
  	\dys  |D_H u_k|^{2}{\rm d}\xi \tows \delta_{\xi_{\rm o}} \ \text{in} \ \mea,
  	$$
  	with~$\delta_{\xi_{\rm o}}$ being the Dirac mass at~$\xi_{\rm o}$.
  \end{theorem}

\vspace{1mm}
\section{Struwe's Global Compactness in the Heisenberg group}\label{sec_glo}
Since the seminal paper~\cite{Str84} by Struwe, the celebrated Global Compactness in the Sobolev space $H^1$ have become a fundamental tool in Analysis which have been proven to be crucial in order to achieve various existence results, as e.~\!g. for ground states solutions for nonlinear Schr\"odinger equations,  for prescribing $Q$-curvature problems, for solutions of Yamabe-type equations in conformal geometry, for harmonic maps from Riemann surfaces into Riemannian manifolds, for Yang-Mills connections over four-manifolds, and many others. The involved literature is really too wide to attempt any reasonable account here. In Theorem~\ref{thm_glob_comp} below, we will state the counterpart of Struwe's Global Compactness in the Heisenberg framework. 

\vspace{2mm}
In order to precisely state such a result, consider  for any fixed~$\lambda\in\r$ the problem,
\begin{equation}\label{plambda}
	-\Delta_Hu-\lambda u-|u|^{2^*-2}u=0\qquad\mbox{in } (\Sc(\Omega))', \tag{$P_\lambda$}
\end{equation}
together with its corresponding Euler--Lagrange energy functional~$\El:\Sc(\Omega)\to\r$ given by
\begin{equation*}\label{Elambda}
	\El(u) =\frac12 \int_{\Om}|D_H u|^2 	\,{\rm d}\xi -\frac{\lambda}{2}\int_{\Om}|u|^2	\,{\rm d}\xi-\frac{1}{2^*}\int_{\Om}|u|^{2^*}	\,{\rm d}\xi.
\end{equation*}
Consider  also the following limiting problem,
\begin{equation}\label{pzero}
	-\Delta_Hu-|u|^{2^*-2}u=0\qquad\mbox{in } (\Sc(\Om_{\rm o}))',\tag{$P_0$}
\end{equation}
where $\Om_{\rm o}$ is either a half-space or the whole~$\h$;
i.~\!e., the Euler-Lagrange equation which corresponds to the energy functional~$\Es: \Sc(\Om_{\rm o})\to\r$,
\begin{equation*}\label{Estar}
	\Es(u)=\frac12 \int_{\Om_{\rm o}}|D_H u|^2 	\,{\rm d}\xi -\frac{1}{2^*}\int_{\Om_{\rm o}}|u|^{2^*}	\,{\rm d}\xi.
\end{equation*}
 
\begin{theorem}[See Theorem~1.3 in~\cite{PPT23}] 
\label{thm_glob_comp} 
	Let~$\{u_k\}\subset \Sc(\Omega)$ be a Palais-Smale sequence for~$\El$; i.~\!e., such that
	\begin{eqnarray*}
		&&\El(u_k)\leq c\quad \mbox{for all }k,\label{PS1}\\
		&&d\El(u_k) \rightarrow 0\quad \mbox{as } k\to\infty \quad \mbox{in }(\Sc(\Omega))'\label{PS2}.
	\end{eqnarray*}
	Then, there exists a (possibly trivial) solution~$\uz\in \Sc(\Omega)$ to~\eqref{plambda} such that, up to a subsequence, we have
	$$
	u_k\rightharpoonup\uz\quad \mbox{as } k\to\infty \quad\mbox{in }\Sc(\Omega).
	$$
	Moreover, either the convergence is strong or there is a finite set of indexes~$\Irm=\{1,\dots,J\}$ such that for all~$j\in\Irm$ there exist a
	nontrivial solution~$\uj\in \Sc(\Om_{\rm o}^{(j)})$ to~\eqref{pzero} with $\Om_{\rm o}^{(j)}$ being either a half-space or the whole~$\h$,
	a sequence of nonnegative numbers~$\{\lambdakj\}$ converging to zero and a sequences of points~$\{\xikj\}\subset\Om$ such that, for a renumbered subsequence, we have for any~$j\in\Irm$
	$$
	\ukj(\cdot):=\lambdakj^{\frac{Q-2}{2}}u_k\big(\tau_{\xikj}\big(\delta_{\lambdakj}(\cdot)\big)\big) \rightharpoonup \uj(\cdot)\quad\mbox{in }\Sc(\h) \quad\mbox{ as }\to\infty.
	$$
	In addition, as~$k\to\infty$ we have
	\begin{eqnarray*}
		&&u_k(\cdot)=\uz(\cdot)+\sum_{j=1}^{J}\lambdakj^{\frac{2-Q}{2}}u_k\big(\delta_{1/\lambdakj}\big(\tau_{\xikj}^{-1}(\cdot)\big)\big)+o(1) \, \mbox{ in }\Sc(\h);\label{propr1}\\
		&&\left|\log{\frac{\lambda_k^{(i)}}{\lambdakj}}\right|+\left|\delta_{1/\lambdakj}\big(\xikj^{-1}\circ \xi_k^{(i)}\big)\right|_{\h}\to\infty\quad\mbox{for }i\neq j,\ \,i,j\in\Irm;\label{propr2}\\
		&&\|u_k\|_{\Sc}^2=\sum_{j=1}^{J}\|\uj\|_{\Sc}^2+o(1); 
		\\
		&&\El(u_k)=\El(\uz)+\sum_{j=1}^{J}\Es(\uj)+o(1) 
	\end{eqnarray*}
\end{theorem}
\mbox{}\vspace{-2mm} \\ In the display above, given~$\xi'\in\h$, we denoted by~$\tau_{\xi'}$ {\it the left translation} defined by
$ \tau_{\xi'}(\xi):=\xi'\circ\xi$ for all $\xi \in \h$.

 \vspace{2mm}
 
The original proof by Struwe in~\cite{Str84} consists in a subtle analysis concerning how the Palais-Smale condition does fail for the functional~$\mathcal{E}^\ast$, based on rescaling arguments, used in an iterated way to extract convergent subsequences with nontrivial limit, together with some slicing and extension procedures on the sequence of approximate solutions to~\eqref{plambda}. Such a proof revealed to be very difficult to extend to different frameworks, and the aforementioned strategy seems even more cumbersome to be adapted to the Heisenberg framework considered here. For this, we completely changed the approach to the problem, and we proved how to deduce the results in Theorem~\ref{thm_glob_comp} in quite a simple way by means of the
so-called {\it Profile Decomposition}, firstly proven by  G\'erard  
 for bounded sequences in the fractional Euclidean space~${H}^s$, and extended to the Heisenberg framework by~\cite{Ben08}. This is in clear accordance with the strategy in~\cite{PP15}; see
 the related result in the fractional Heisenberg framework in~\cite{GMM18}.
\begin{remark}
The limiting domain~$\Om_{\rm o}$ in   {\rm Theorem~\ref{thm_glob_comp}}  can be either the whole~$\h$ or a half-space. On the contrary, in the original proof in the Euclidean case by~Struwe~{\rm (\cite{Str84})} one can exclude the existence of nontrivial solutions to the limiting problem in the half-space by   Unique Continuation and Pohozaev's Identity. Such a possibility can not be a priori excluded  in the sub-Riemannian setting, even in the very special case when a complete characterization of the limiting set is possible under further regularity assumptions on~$\Omega$. Indeed, in the Heisenberg framework, a very few nonexistence results are known, basically only in the case when the domain reduces to a~half-plane parallel or perpendicular to the group center; see~\cite{CU01}. We also refer to the last paragraphs in~\cite[Section 5] %ALERT
{PPT23} for further details.
\end{remark}

\vspace{1mm}
\section{Asymptotics of the optimal functions}\label{sec_asy}
We present an asymptotic control of the maximizing sequence~$u_\eps$ for~$\Ssub$ in~\eqref{sobolev} via the Jerison\,\&\,Lee extremals. This is shown in Theorem~\ref{han} below,  
 which will be one of the key in the proof of the localization of the concentration result presented in Section~\ref{sec_con} below and it could be also
  useful to investigate further properties related to subcritical Folland-Stein-Sobolev embeddings. 
    \begin{theorem}[See Theorem~1.2 in~\cite{PP23h}]\label{han}
   	Let~$\Omega \subset \h$ be a  
	smooth 
	bounded  domain such that  
    \begin{equation*}\label{non_deg_Om}
     \liminf_{\rho \to 0}\frac{\snr{(\h \setminus \Om) \cap B_\rho(\xi)}}{\snr{B_\rho(\xi)}}  >0 \ \, \forall \xi \in \partial \Om.
    \end{equation*}
    Then, for    each~$0<\eps< 2^*-2$ letting~$u_\eps \in \Sc(\Om)$  being a maximizer for~$\Sob_\eps$, there exist~$\{\eta_\eps\} \subset \Om$,~$\{\lambda_\eps\} \subset \r^+$ such that, up to choosing~$\eps$ sufficiently small, we have that
   	\begin{equation*}\label{bound_max_seq}
   		u_\eps \lesssim  \,U_{\lambda_\eps,\eta_\eps} \qquad ~\textrm{on}~\Om,
   	\end{equation*}
   	where~$U_{\lambda_\eps,\eta_\eps}=U \left(\delta_{1/\lambda_\varepsilon}\big(\tau_{\eta_\eps}(\xi)\big)\right)$ are the {\rm Jerison{\,\rm\&\,}Lee} extremal functions, 
	 and the sequences~$\{\eta_\eps\}$ and~$\{\lambda_\eps\}$ satisfy
   	\begin{equation*}\label{conv_xi_eps}
   	\eta_\eps \sim \,  \xi_{\rm o} \quad \textrm{and} \quad
   		\lambda_\eps^\eps \sim 1 \quad \textrm{as} \ \eps \searrow 0,
   	\end{equation*}  
    with~$\xi_{\rm o}$ being the concentration point given in {\rm Theorem~\ref{cor_concentration}}.
   \end{theorem}
The result in Theorem~\ref{han} above reminds 
to the literature following the pioneering work in the Euclidean framework due to Aubin and Talenti, and in such a framework it is fundamental in the proof of a precise conjecture about the localization of the concentration point~$\xi_{\rm o}$ given in Corollary~\ref{cor_concentration} by~Han in~\cite{Han91}.  
In the proof of Theorem~\ref{han} in~\cite{PP23h} in the sub-Riemannian framework we are dealing with, one has also to deal with the fact that, in strong contrast with the Euclidean setting, the Jerison\,\&\,Lee extremals cannot be reduced to functions depending only on the standard Kor\'anyi gauge. For this, such a proof will require a delicate strategy which makes use and refines the concentration result obtained via the $\Gamma$-convergence result in Theorem~\ref{the_gamma-intro} in order to detect the right involved scalings~$\eta_\eps$ and~$\lambda_\eps$. Also the Global Compactness-type result presented in Section~\ref{sec_glo} is  needed.

\vspace{1mm}
\section{Localization of the energy concentration}\label{sec_con}
A natural question arises: can 
the blowing up be localized; i.~\!e.,  is the concentration point~$\xi_{\rm o}$ in Theorem~\ref{cor_concentration} in~Section~\ref{sec_gamma}
 related in a specific way to the geometry of the domain~$\Omega$\,?
 \\ In the Euclidean framework, under standard regularity assumptions,  
  Han (\cite{Han91}) and Rey~(\cite{Rey89}) proved the connection with the Green function associated to the domain~$\Om$ by answering to a famous conjecture by Brezis and Peletier~(\cite{BP89}), who had previously investigated the spherical domains setting. The involved proofs strongly rely on the regularity of Euclidean domains, which is in clear contrast with the complexity of the underlying sub-Riemannian geometry here; as well-known, even if the domain~$\Omega$ is smooth, the situation is drastically different because of the possible presence of characteristic points on the boundary~$\partial\Omega$. From one side, near those characteristic points -- as firstly discovered by Jerison  -- even harmonic functions on the Heisenberg group can encounter a sudden loss of regularity;
from the other side, one did not want to work in the restricted class of domains not having characteristic points. 
  In order to deal with those specific difficulties, it is thus quite natural to work under the assumption that the domain~$\Om$ is {\it geometrical regular near its characteristic set} as given by Definition~\ref{def_nnta} below. 
In forthcoming Theorem~\ref{thm_green} we state  the expected localization result for the concentration point~$\xi_{\rm o}$~of the maximizing sequence~$u_\eps$ in terms of the Green function associated with the domain~$\Om$, in turn establishing the validity of the aforementioned Brezis-Peletier conjecture in the Heisenberg group.
 \vspace{1mm}
 
As customary,  denote by~$\mathcal{D}$ the infinitesimal generator of the one-parameter group of non-isotropic dilations~$\{\delta_\lambda\}_{\lambda>0}$ in~\eqref{def_philambda}; that is,
  \begin{equation}\label{dilation_vector_field}
  	\mathcal{D} := \sum_{j=1}^n \big(x_j\partial_{x_j}+y_j \partial_{y_j}\big) +2t\partial_t.
  \end{equation}
  
   \begin{definition}[\bf $\delta_\lambda$-starlike sets]
  	Let~$\Om$ be a~$C^{1}$~connected open set of\,~$\h$ containing the group identity~$\mathfrak{e}$. We say that $\Om$ is \textup{$\delta_\lambda$-starlike} {\rm (}with respect to the identity~$\mathfrak{e}${\rm)} along a subset~$K \subseteq \partial \Om$ if
  	\[
  	\langle \mathcal{D},\mathfrak{n}\rangle (\eta) \ge 0, 
  	\]
  	at every~$\eta \in K$;  in the display above~$\mathfrak{n}$ indicates the exterior unit normal to~$\partial \Om$. 
  	
  	We say that~$\Om$ is \textup{uniformly $\delta_\lambda$-starlike} {\rm (}with respect to the identity~$\mathfrak{e}${\rm)} along~$K$ if there exists~$\alpha_\Om >0$ such that, at every~$\eta \in K$,
  	\[
  	\langle \mathcal{D},\mathfrak{n}\rangle (\eta) \ge \alpha_\Om.
  	\]
  	A domain as above~$\Om$ is {$\delta_\lambda$-starlike} (uniformly {$\delta_\lambda$-starlike}, respectively) with respect to one of its point~$\zeta \in  \Om$ along~$K$ if~$\tau_{\zeta^{-1}}(\Om)$ is {$\delta_\lambda$-starlike} (uniformly {$\delta_\lambda$-starlike}, respectively) with respect to the origin along~$\tau_{\zeta^{-1}}(K)$.
  \end{definition}

  Given a domain~$\Om\subset\h$, we recall that its {\it characteristic set}\,~$\Car$, the collection of all its characteristic point, is given by
  \[
  \Car :=\Big\{\xi \in \partial \Om \, | \, Z_j (\xi) \in T_\xi (\partial\Om), \, \textrm{for}~j=1,\dots,2n\Big\}.
  \]

  We now recall the definition of  regular domains in accordance with the by-now classical paper~\cite{GV00}.   
  \begin{definition}[See Definition~2.2 in~\cite{PP23h}] \label{def_nnta}
%ALERT

    	A smooth domain~$\Om\subset\h$ such that~$\partial \Omega$ is an orientable hypersurface is {\rm  
    		``geometrical regular near its characteristic set''} if
    	the following conditions hold true,
    	\begin{itemize}
    		\item[($\Om1$)]
    		There exist~$\varPhi \in C^\infty(\h)$,~$c_\Om >0 $ and~$\rho_\Om \in \r$ such that
    		\[
    		\Om := \big\{\varPhi < \rho_\Om\big\}, \quad \textrm{and} \quad \snr{D \varPhi} \geq c_\Om.
    		\]
    		\vspace{1mm}
    		\item[($\Om2$)]
    		For any~$\xi \in \partial \Om$ it holds
    		\[
    		\liminf_{\rho\to0^+}\frac{\snr{(\h \smallsetminus \Om) \cap B_\rho(\xi)}}{ \snr{B_\rho(\xi)}}>0.
    		\]
    		    		\vspace{1mm}
    		\item[($\Om3$)] 
    		There exist~$M_\Om$ such that
    		\[
    		\Delta_H \varPhi \ge \frac{4\snr{z}}{M_\Om } \langle D_H \varPhi, D_H \snr{z} \rangle \quad \textrm{in}~\omega,
    		\]
    		where~$\omega$ is an interior neighborhood of~$\Car$.
    		    		\vspace{1mm}
    		\item[($\Om4$)]
    		$\Om$~ is {$\delta_\lambda$-starlike} with respect to one of its point~$\zeta_{\rm o} \in \Om$ and uniformly {$\delta_\lambda$-starlike} with respect to~$\zeta_{\rm o}$ along~$\Car$.
    	\end{itemize}
\end{definition}

 We are finally in the position to state the localization result.%ALERT
    \begin{theorem}[See Theorem~1.3 in~\cite{PP23h}]\label{thm_green}
    Consider a bounded domain~$\Om\subset\h$  geometrical regular near its characteristic set,  
     and  
    let~$u_\eps \in \Sc(\Om)$  be a maximizer for~$\Sob_\eps$. Then, up to subsequences, $u_\eps$~concentrates at some point~$\xi_{\rm o} \in  \Om$ such that
       \begin{equation}\label{robin_condition}
      	\int_{\partial \Om} \snr{D_H G_\Om(\cdot,\xi_{\rm o})}^2 \langle\mathcal{D},\mathfrak{n}\rangle \, {\rm d}\mathscr{H}^{Q-2}=0,
      \end{equation}
      with~$G_\Om(\cdot;\xi_{\rm o})$ being the Green function associated to~$\Om$ with pole in~$\xi_{\rm o}$, and $\mathcal{D}$ being the infinitesimal generator of the one-parameter group of non-isotropic dilations in the Heisenberg group defined in~\eqref{dilation_vector_field}. 
    \end{theorem}
    
    %ALERT
    The proof  can be found in~Section~7 in~\cite{PP23h}; it involves all the results stated in the preceding sections together with other general tools in the sub-Riemaniann framework, as e.\!~g., maximum principles, Caccioppoli-type estimates, $H$-Kelvin transform, boundary Schauder-type regularity estimates, as well as with a fine  boundary  analysis of the solutions to subcritical Yamabe equations.
    We refer also to the interesting related result in~\cite{MMP13} in the case of domains with no characteristic points.

\section*{Acknowledgements}
The authors are member of Gruppo Nazionale per l'Analisi Matematica, la Probabilit\`a e le loro Applicazioni (GNAMPA) of Istituto Nazionale di Alta Matematica ``F.~Severi'' (INdAM), whose support is acknowledged.  The authors are also supported by INdAM Project ``Fenomeni non locali in problemi locali", \,CUP\_\,E55F220\break 00270001.
The second author is also supported by the Project ``Local vs Nonlocal: mixed type operators and nonuniform ellipticity", \!CUP\_\,D91B21005370003.

\end{document}